\documentclass[journal,draftclsnofoot,onecolumn]{IEEEtran}
\usepackage{amsmath}
\usepackage{amssymb}
\usepackage{amsthm}
\usepackage[shortlabels]{enumitem}
\usepackage{graphicx}
\usepackage{url}
\usepackage{color}
\usepackage{cite}
\usepackage{mathtools}
\usepackage{tikz}
\usetikzlibrary{shapes.misc,decorations.pathreplacing,decorations.markings,patterns,arrows.meta}

\newcommand{\RR}{{\mathbb R}}

\newcommand{\cB}{{\mathcal{B}}}
\newcommand{\cC}{{\mathcal{C}}}
\newcommand{\cM}{{\mathcal{M}}}

\newcommand{\Prob}{{\mathbb{P}}}

\newcommand{\cP}{{\mathcal{P}}}

\newcommand{\cF}{{\mathcal{F}}}

\newcommand{\cL}{{\mathcal{L}}}

\newcommand{\cT}{{\mathcal{T}}}
\newcommand{\argmax}{\mathop{\mathrm{argmax}}}

\tikzset{
    partial ellipse/.style args={#1:#2:#3}{
        insert path={+ (#1:#3) arc (#1:#2:#3)}
    }
}

\newtheorem{theorem}{Theorem}
\newtheorem{proposition}[theorem]{Proposition}
\newtheorem{remark}[theorem]{Remark}

\newtheorem{example}[theorem]{Example}
\newtheorem{problem}[theorem]{Problem}

\newtheorem{definition}[theorem]{Definition}

\author{Melkior Ornik and Ufuk Topcu
\thanks{M.~Ornik is with the Institute for Computational Engineering and Sciences, University of Texas at Austin. e-mail: mornik@ices.utexas.edu}%
\thanks{U.~Topcu is with the Department of Aerospace Engineering and Engineering Mechanics and the Institute for Computational Engineering and Sciences, University of Texas at Austin. e-mail: utopcu@utexas.edu}%
}
\title{Deception in Optimal Control}
\date{}

\begin{document}
\maketitle

\begin{abstract}
In this paper, we consider an adversarial scenario where one agent seeks to achieve an objective and its adversary seeks to learn the agent's intentions and prevent the agent from achieving its objective. The agent has an incentive to try to deceive the adversary about its intentions, while at the same time working to achieve its objective. The primary contribution of this paper is to introduce a mathematically rigorous framework for the notion of deception within the context of optimal control. The central notion introduced in the paper is that of a belief-induced reward: a reward dependent not only on the agent's state and action, but also adversary's beliefs. Design of an optimal deceptive strategy then becomes a question of optimal control design on the product of the agent's state space and the adversary's belief space. The proposed framework allows for deception to be defined in an arbitrary control system endowed with a reward function, as well as with additional specifications limiting the agent's control policy. In addition to defining deception, we discuss design of optimally deceptive strategies under uncertainties in agent's knowledge about the adversary's learning process. In the latter part of the paper, we focus on a setting where the agent's behavior is governed by a Markov decision process, and show that the design of optimally deceptive strategies under lack of knowledge about the adversary naturally reduces to previously discussed problems in control design on partially observable or uncertain Markov decision processes. Finally, we present two examples of deceptive strategies: a ``cops and robbers'' scenario and an example where an agent may use camouflage while moving. We show that optimally deceptive strategies in such examples follow the intuitive idea of how to deceive an adversary in the above settings.

\end{abstract}

\section{Introduction}

The concept of deception is naturally present in a variety of contexts that have an adversarial element. Examples include cybersecurity \cite{CarGro11, PanLiu12, Kwoetal13}, bio-inspired robotics \cite{ShiArk14}, genetic algorithms \cite{Whi91}, vehicle decision making \cite{McESin05}, warfare strategy --- a particularly voluminous study of the role of deception in war has been made in \cite{Wha69} --- and interpersonal relationships \cite{Met89}. A {\em deceptive strategy} employed by an agent in an adversarial setting rests on dual goals of the agent:
\begin{enumerate}[1)]
\item achieving its objective,
\item modifying the adversary's beliefs about the nature of that objective --- for instance, objective location, distance, or reward attained at an objective.
\end{enumerate}
The desire to modify the adversary's beliefs is motivated by the assumption that the adversary would be able to block, or modify, access to the objective if it correctly identified its nature. Thus, the success of an agent at achieving its control objective often depends on the agent's ability to hide its true intentions from its adversary while still proceeding towards its objective, and satisfying any other constraints that may be placed on its behavior (e.g., safety constraints).

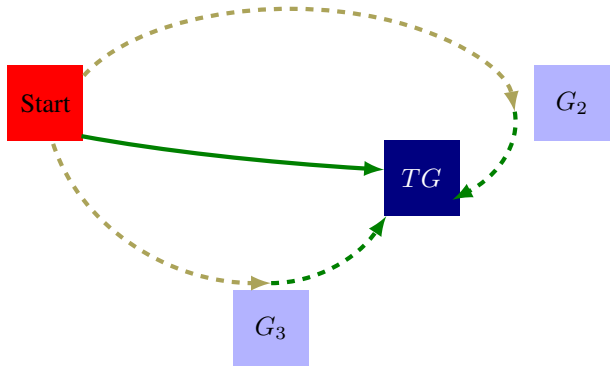
\begin{figure}[ht]
\centering
\begin{tikzpicture}
\fill[blue!50!black] (6,4) rectangle (5,5);
\fill[blue, opacity=0.3] (3,2) rectangle (4,3);
\fill[blue, opacity=0.3] (8,5) rectangle (7,6);
\fill[red] (0,5) rectangle (1,6);
\node at (5.5,4.5) {\color{white} $TG$};
\node at (3.5,2.5) {$G_3$};
\node at (7.5,5.5) {$G_2$};
\node at (0.5,5.5) {Start};

\draw[ultra thick, yellow!60!black, -latex, dashed] (3.75,5.25) [partial ellipse=156:5:3cm and 1.5cm];
\draw[ultra thick, green!50!black, -latex, dashed] (3.75,5.25) [partial ellipse=5:-44:3cm and 1.5cm];
\draw[ultra thick, green!50!black, -latex] (8.75,6) [partial ellipse=219:248:10cm and 1.5cm];
\draw[ultra thick, yellow!60!black, -latex, dashed] (3.5,5.6) [partial ellipse=195:270:3cm and 2.5cm];
\draw[ultra thick, green!50!black, -latex, dashed] (3.5,5.6) [partial ellipse=270:320:2cm and 2.5cm];
\end{tikzpicture}
\caption{An illustration of the running example of this paper. The agent seeks to take a path from the starting area, denoted in red, to its true goal ($TG$), denoted in dark blue. An adversary seeks to determine which of the three possible blue goals ($TG$, $G_2$, $G_3$) is the agent's true goal. If the agent takes the direct path (denoted in solid green) from the start to the true goal, the adversary can easily correctly infer its intentions. If the agent, however, first begins moving towards one of the false goals (dashed yellow paths) the adversary may assume that this goal is the agent's objective, and may only have very little time to change its opinion, after the agent does not end up going into a false goal, and turns towards the true goal. Thus, the agent's top and bottom paths are deceptive.}
\label{fig0}
\end{figure}

A simple example for a deceptive setting, which will serve as the running example throughout this paper, is given in Figure \ref{fig0} --- we refer to it as {\em cops and deceptive robbers}. In the setting illustrated by Figure \ref{fig0}, the agent (``robbers'') seeks to move to a particular area of a state space (``bank''), which holds a reward. An adversary (``cops'') knows that the agent is seeking to reach one of several possible objectives (``banks''), but does not know which one. By observing the agent's movement, the adversary seeks to learn the agent's intentions, and change the nature of the objective (``set a trap''). Hence, it is in agent's interest to make the adversary's beliefs about the agent's intentions as incorrect as possible, while still ultimately reaching its objective.

Within the context of this paper, we are primarily interested in encoding the agent's strategy as controlled movements on an state space over a period of time. Such an approach allows us to naturally discuss scenarios such as our running example, where the agent is moving across physical space, but also explore more general settings where agent changes its behavior, values, and objectives over time. We encode the agent's objectives through a reward function that should be maximized. Additionally, we allow for existence of possible a priori constraints on the agent's behavior, notably temporal logic specifications. 

Despite aforementioned wide interest in deception and design of deceptive strategies, mathematically formal definition of deception in existing literature has largely been limited to application-specific settings. For example, \cite{PanLiu12} considers deception solely in the context of cyber-attacks on networked systems, \cite{Yav87} deals with pursuit-evasion scenarios, while \cite{Casetal04} discusses area denial. Significant theoretical literature on deception exists in the context of abstract two-agent games; we turn the reader's attention to \cite{Hesetal00,WagArk09,Arketal12}. Nonetheless, the considered scenarios often do not  account for time-varying behavior. Instead, the relationship between the agent and its adversary are constrained to single-stage games, where the level of trust between the parties do not evolve over time. Multi-stage two-agent games are considered in \cite{Pau15,Sin06,EttJeh10}. In particular, \cite{Pau15} primarily concentrates on a framework where the adversary's learning period (in which the adversary may be deceived) is separated from the ``main'' period during which the movement takes place, and is thus significantly different from the framework that we aim to develop, which involves an agent deceiving the adversary {\em while} attempting to achieve its objective. The frameworks developed in \cite{Sin06} and \cite{EttJeh10} bear more similarity with ours. Nonetheless, they differ in substantial elements. First of all, both papers use two agents which are both able to take actions, and have separate objectives, during the system run. In our setting, the deceiver (i.e., agent) is the only one who can take actions, while the adversary solely observes, learns from a predetermined model, and influences the agent's rewards, but does not itself move in the state space. Additionally, \cite{Sin06} frames deception as exploitation of the adversary's lack of knowledge about the system state. Our intuition is different. We allow for the possibility that the adversary is able to see the agent's states and actions at all times. The unknown element are the agent's intentions. In our technical approach to deception, we make the adversary's belief of these intentions a part of the enlarged state space. Thus, the framework of \cite{Sin06} could perhaps be related to ours, if the state space in \cite{Sin06} is understood as this enlarged state space. Nonetheless, even in the enlarged state space, the agent and the adversary do not both move in the same way: the agent's movements might change the adversary's beliefs, but the adversary performs no movements on its own.

Paper \cite{EttJeh10}, on the other hand, primarily focuses on two agents that are effectively symmetric in the sense that they attempt to deceive each other --- as opposed to our setup where one agent deceives and the other is being deceived. It then makes significant assumptions on the behavior, i.e., rationality, of both agents. It then shows that there exists a mutually optimal deception strategy --- a partial information version of a Nash equilibrium.  In contrast, we are not interested in an equilibrium setting, instead wanting to concentrate solely on a single agent's choices at deceiving its adversary. Additionally, \cite{EttJeh10} deals solely with discrete-time games, and primarily does so in the context of agents proceeding through multiple rounds of negotiations towards an outcome (e.g., trade agreement). Naturally, there is no explicit notion of a state space for either agent to move in, and no opportunity to easily encode specifications, given in, e.g., temporal logic, that limit the agent's control policies. While, as with \cite{Sin06}, it may be possible to partly reframe the framework of this paper into the setting of \cite{EttJeh10}, and vice versa, the setting of \cite{EttJeh10} is not naturally amenable to scenarios involving an agent moving across a state space in order to collect rewards encoded in such a space.

The objective of this paper is to formalize deception, and deceptive strategies, within the framework of optimal control, as well as discuss optimal design of deceptive strategies for a wide class of scenarios. Since a critical component of deception is modifying adversary's opinions, we begin our discussion by formalizing the notion of adversary's {\em belief space} and {\em belief-induced rewards}. We use these notions to define deception and optimal deceptive strategies for general control systems in Section \ref{defdec}. As it is natural that the deceiving agent might not know everything about adversary's beliefs, Section \ref{lack1} briefly discusses several models of lack of agent's knowledge about the adversary. Concretizing the general framework of Section \ref{defdec}, Section \ref{mdpsec} concentrates on deceptive strategies on Markov decision processes, and Section \ref{lack2} considers design of optimal deceptive strategies for the previously identified classes of lack of knowledge about the adversary. In the latter part of the paper, we examine two particular examples of deception: the running example of a ``cops and robbers'' setting in which the adversary knows multiple possible candidate objectives, but does not know which is the true objective, is discussed in Section \ref{par}, and a setting in which the adversary is aware of the objective, but does not possess complete information about the agent's position, is discussed in Section \ref{camor}. In Section \ref{conc} we briefly describe future work concerning deceptive scenarios, design of optimal deceptive strategies, and models of the agent's lack of knowledge. 

\section{Definition of Deception}
\label{defdec}

We now define deception and deceptive strategies in control systems with reward-based objectives. In order to motivate our definition of deception in an adversarial scenario where the agent has an objective, and an adversary is attempting to learn the agent's intentions and influence its achievement of the objective, we first consider a simplified scenario with the same control objective, but without an adversary. We refer to such a scenario as {\em nominal}. 

Let $\cC$ be a system evolving on a set of times $\cT$, and describing a single agent's behavior as it moves across a state space $S$ using controls from the action set $A$. For the purposes of this section, we make no assumptions on the structure of $\cC$. Solely for notational reasons, we consider $\cT$ to be a discrete set, but we emphasize that all the notions can be equivalently posed when $\cT$ is an interval.
Assume that system $\cC$ comes equipped with a {\em nominal reward} function $R:S\times A\times \cT\to\RR$. The control objective of an agent evolving in $\cC$ is to maximize its accumulated reward over some, potentially infinite, period of time.

\begin{definition}
\label{def1}
The {\em nominal optimal control policy} is given by
\begin{equation}
\label{nomeq}
\argmax_{a_t} \sum_{t\in \cT} R(s_t,a_t,t)\textrm{,}
\end{equation} where $s_t$ is the agent's state at time $t$ and $a_t$ is the agent's action at that time. If the dynamics on $\cC$ are stochastic, as, for example, in the case of Markov decision processes (MDPs) discussed in subsequent sections, \eqref{nomeq} is usually replaced by
\begin{equation}
\label{nomeqex}
\argmax_{a_t} \mathbb{E}\left[\sum_{t\in \cT} R(s_t,a_t,t)\right]\textrm{,}
\end{equation}
given that the actual sum is impossible to calculate in advance.
\end{definition}

We note that, in Definition \ref{def1}, it is possible to add constraints on the choice of actions available to the agent at any given time $t$. Namely, instead of allowing that $a_t$ be any element of $A$ at any time $t$, we may require $a_t\in A_t\subseteq A$, where $A_t$ may depend on the agent's history, i.e., on previous actions $a_{t'}$ for $t'\leq t$, previous states $s_{t'}$, as well as $s_t$. For instance, in an MDP framework discussed in subsequent sections, such limitations naturally arise out of temporal logic specifications on the agent's path. We briefly discuss such constraints, and their implications to decepvie strategies, in the following section. However, in terms of Definition \ref{def1} and the subsequent formal discussion in this section, the only change that such specifications impose is that the arguments in $\argmax$ satisfy $a_t\in A_t$ for all $t\in\cT$. Thus, we omit future reference to such specifications in the present section.



Having defined a control objective and a nominal optimal control policy for an agent without the presence of an adversary, we now seek to formalize the adversary's role. In our framework, the adversary has two salient properties:
\begin{enumerate}[1)]
\item belief about the agent's intentions, which may change over time, and
\item influence of the adversary's belief on the agent's actual collected reward.
\end{enumerate}

We note that we have not concretized the meaning of a {\em belief} in 1), nor the meaning of an {\em intention}. Informally, we consider intention to be any property of the agent, or agent's policy, that is important to the adversary. For instance, it may be the agent's objective, agent's next action, or the agent's accumulated reward. A belief is, then, an assertion on the set of agent's intentions.

Formally, we simply define a {\em belief} $B_t$ at time $t\in\cT$ as an element of some domain $\cB$, which we refer to as the {\em belief space}. As stated above, $\cB$ can generally be any set that in some way describes the elements of agent's behavior that are important to the adversary. A concrete instantiation of $\cB$ depends on the exact setting that we are dealing with. For instance, in the context of our running example, one possible model, which we revisit later in the text, is that $\cB$ consists of all states that are the agent's possible goals. Another model would define $\cB$ as the set of probability distributions on the set of possible objectives, thus allowing that the adversary is uncertain about the agent's objective.

Having formally described the adversary's beliefs, we now move to property 2), i.e., describe how those beliefs change the reward associated system $\cC$. In our running example, the knowledge of the cops about the robbers' intentions will change the payoff that the robbers' get for reaching the objective: instead of robbing the bank, they will be caught. In other words, instead of collecting a reward $r>0$, the agent will collect a reward $r'<0$. We generalize this notion by introducing {\em belief-induced rewards}.

\begin{definition}
\label{def2}
A {\em belief-induced reward} function is a map given by $L:S\times\cB\times A\times \cT\to\RR$. 
\end{definition}

While there is no requirement that the nominal reward $R$ and belief-induced reward $L$ from Definition \ref{def2} are in any way related, the motivating setting would imply that $L$ is in some way a modification of $R$. For instance, one could consider $L(s,B,a,t)=R(s,a,t)+H(B)$, where $H$ is a function that depends solely on the adversary's belief. However, we will not be a priori assuming any formal relationship between $L$ and $R$ in the theoretical results of this paper. Instead, we will be dealing with optimal design of agent's policy with respect to a general reward function $L$, while making use of the relationship between $R$ and $L$ when analyzing the effectiveness of deception for scenarios explored in Section \ref{par} and Section \ref{camor}.

Analogously to the nominal optimal control policy from Definition \ref{def1}, a belief-induced reward from Definition \ref{def2} yields an optimal belief-induced control policy.

\begin{definition}
\label{def2b}
The {\em optimal belief-induced control policy} is given by
\begin{equation}
\label{acteq}
\argmax_{a_t} \sum_{t\in \cT} L(s_t,B_t,a_t,t)\textrm{,}
\end{equation} where $s_t$ is the agent's state at time $t$, $a_t$ is the agent's action at that time, $B_t$ is the adversary's belief at time $t$.
\end{definition}

Our running example can now finally be more formally defined as follows. In it, we assume that the agent has a single true goal $TG\in S$. The adversary possesses a set of possible goals, denoted by $\{TG,G_2,G_3\}$ in Figure \ref{fig0}, but does not know which one is the true goal. If the adversary's belief of the agent's true goal is incorrect, i.e., the agent successfully fooled the adversary, then the agent collects a positive reward for reaching the true goal. On the other hand, if the adversary's belief is correct, then the agent collects a negative reward for reaching the true goal. The above setup is formalized as follows.

\begin{example}[Cops and deceptive robbers]
\label{exa0}
Let $\cC$ be a control system with a corresponding state space $S$, action set $A$, and set of times $\cT$. Let $TG\in S$ be the agent's true goal, and $\{G_1,\ldots,G_k\}\subseteq S$, with $TG\in\{G_1,\ldots,G_k\}$, be the set of states that the adversary believes are the possible agent's objectives. Then, define $\cB=\{G_1,G_2,\ldots,G_k\}$. Finally, let the belief-induced reward $L:S\times \cB\times A\times \cT\to\RR$ be defined by 
\begin{equation*}
L(s,B,a,t)=\begin{cases}
0 & \textrm{for all } s\in S\backslash\{TG\}, B\in\cB, a\in A, t\in T\textrm{,} \\
L^+ & \textrm{for all } s=TG, B\in\cB\backslash\{TG\}, a\in A, t\in\cT\textrm{.} \\
L^- & \textrm{for all } s=TG, B=TG, a\in A, t\in\cT\textrm{,}
\end{cases}
\end{equation*}
where $L^+$ and $L^-$ are real numbers with $L^+>0$ and $L^-<0$. 
The optimal belief-induced policy for the robbers is then given by \eqref{acteq}.
\end{example}

\subsection{Optimal Control Design}
\label{optcd}

Having established the notion of an optimal belief-induced strategy, we now consider the problem of determining an optimal belief-induced policy for the agent. In the setting introduced above we generally assume that the adversary is not controlled; the control system $\cC$ was defined solely for the agent, while the adversary's beliefs change according to a predetermined learning mechanism. While the adversary's belief evolves on the belief space $\cB$, the adversary has no goal that it attempts to reach. Its movement is entirely defined by the agent's actions. Such an assumption makes the setting that we are considering less general than the framework of two-player games. On the other hand, having complete knowledge of the adversary's reaction allows us to describe an optimal policy for the agent. We are also able to consider scenarios where the agent has only partial information about the adversary, in which case we can discuss an optimal policies {\em given} the agent's knowledge. We describe such a setting in Section \ref{lack1} and Section \ref{lack2}.

Assuming that the adversary's belief changes over time depending on the agent's trajectory and actions, the control system $\cC$ along with the dynamics on $\cB$ defines a derived belief-induced control system $\cC_\cB$ on $S_\cB=S\times\cB$. Thus, the problem of finding an optimal belief-induced policy described in Definition \ref{def2b} can be understood as a reward maximization problem in a control system $\cC_\cB$ on state space $S_\cB$. Thus, at every time, the agent may commit an action with the purposeful desire to modify the adversary's belief in a way that will lead to an increase in the agent's actual attained reward. Such a decision would necessarily be based on some knowledge of the adversary's current belief, or its belief dynamics, and an optimal belief-induced policy is thus an {\em optimal deceptive policy}. In the context of this paper, we define {\em deception} as any such exploitation of prior or side information that the agent may have on $L$ and dynamics of $\cB$ to better design its control policy. For instance, in the context of our running example, an agent uses deception if it does not merely go towards its true goal in the most direct path, but in some way exploits the knowledge that an adversary is attempting to learn its goal, and may reduce the agent's reward if it learns the goal correctly.

Without additional assumptions on system dynamics in $\cC_\cB$, the problem of finding an optimal policy \eqref{acteq} is an optimal control problem on a product space $S\times\cB$. Solvability of such a problem, or computational difficulty of finding a solution, depends on the details of the system dynamics. For a particularly extensive treatise of optimal control strategies, see \cite{Ber12,Ber17}. We will concentrate on particular models for the state space $S$ and belief space $\cB$ in subsequent sections. However, before moving to that segment of the paper, let us briefly discuss the role of {\em limited knowledge} about the system $\cC_\cB$ in deception.

\subsection{Lack of Knowledge}
\label{lack1}

In Definition \ref{def2b}, we defined an optimal belief-induced, or {\em deceptive}, control policy, and we showed above that computing such a policy is equivalent to a solution of a general optimal control problem in an appropriate domain. However, finding an optimal control policy depends on entirely knowing the dynamics of the belief-induced system $\cC_\cB$, as well as perfectly observing the full system state $(s_t,B_t)$ at all times. Such an assumption is not realistic in many adversarial scenarios --- the adversary may have interest in not divulging its current belief to the agent. Thus, it may be impossible for the agent to devise an objectively optimal deceptive policy; instead, the goal is to devise a deceptive policy that is optimal given the agent's imperfect knowledge.

We consider three categories in which the agent may lack precise knowledge on the belief-induced control system $\cC_\cB$. In particular, the agent might not know:

\begin{enumerate}[1)]
\item \textbf{What the adversary thinks} --- knowledge of the adversary's belief $B_t$ at every time step.
\item \textbf{How the adversary thinks} --- knowledge of the dynamics underlying $B_t$.
\item \textbf{What the adversary does} --- knowledge on how the actual reward $L$ depends on the adversary's beliefs $B$.
\end{enumerate}
This list of categories is not exhaustive, and we briefly discuss some other variants of lack of knowledge in Section~\ref{conc}.

While the definition of deception established in Section \ref{defdec} is valid for any general control system $\cC_\cB$, and, as mentioned, problems involving deception can be generally positioned within the umbrella of optimal control, the above definition is clearly too broad to allow for theoretical and algorithmic results on performance of deceptive agents and design of deceptive strategies. Thus, as mentioned previously, in the following section we make our setting more specific by focusing on Markov decision processes.

\section{Deception in Markov Decision Processes}
\label{mdpsec}

Section \ref{defdec} defined a notion of deception for agents operating on any general state space $S$ endowed with a reward function. As we showed in Section \ref{optcd}, the problem of finding an optimal deceptive policy reduces to a solution of an optimal control problem on a set $S\times\cB$. However, strategies for design of optimal control policies in general state spaces generally significantly depend on the structure of the state space that the agent evolves on. Thus, in the remainder of this paper we assume a particular structure of the state space $S$; namely, we assume that the evolution of an agent is given by a finite-state discrete-time Markov decision process.

A {\em Markov decision process (MDP)} is defined by $\cM=(S,A,P)$, where the state space $A$ and action set $A$ are finite, and the agent's dynamics on $S$ are given by \begin{equation*}
\Prob(s_{t+1}=s')=P(s_t,a,s')\textrm{,}
\end{equation*} where $P:S\times A\times S\to[0,1]$ satisfies $P(s,a,s')\geq 0$ for all $s,s'\in S$, $a\in A$, and $\sum_{s'\in S} P(s,a,s')=1$ for all $s\in S$, $a\in A$. In the remainder of the paper, we will be interested in finite-time system runs, i.e., we assume that the interval $\cT$ of interest is given by $\cT=\{0,1,\ldots,T\}$. We also assume that all elements of $\cM=(S,A,P)$, i.e., state and action sets, as well as transition probabilities, are known to both agent and the adversary.

\subsection{Optimal Deceptive Policy}

As described in Definition \ref{def1}, the agent's nominal objective \eqref{nomeqex} is to maximize 
\begin{equation*}
\mathbb{E}\left[\sum_{t=0}^T R(s_t,a_t)\right]\textrm{,}
\end{equation*}
where $s_t$ is the agent's position at time $t$, and $a_t$ the action that the agent took at time $t$.

\begin{remark}
We note that, in the previous section, $R$ was allowed to depend on time along with $S\times A$. We consider time-invariant rewards in order to parallel standard MDP setups (see, e.g., \cite{Put05} for a detailed introduction to MDP rewards). Nonetheless, time dependence could be easily encoded by considering the state space $\tilde{S}=S\times\cT$.
\end{remark}


As described in Definition \ref{def2}, the agent's nominal reward $R$ is modified by the adversary's beliefs into a function $L:S\times\cB\times A\to\RR$. In the remainder of this paper, we assume that the belief set $\cB$ is finite. Such an assumption ensures that the state space $S\times\cB$ of the belief-induced system $\cC_\cB$ is finite. We emphasize that the theoretical framework that we are developing does not prescribe any further structure of the set $\cB$, i.e., the meaning of the beliefs.

\begin{remark}
If $\cB$ is infinite, the problem of determining an optimal deceptive policy \eqref{acteq} is a problem of optimal control on an infinite state space $S\times\cB$. We point a reader interested in optimal control on infinite-state MDPs to \cite{Sch93}. An infinite-state framework may appear naturally, for instance, if $\cB$ is the set of all possible reward functions $R:S\times A\to\RR$. In that particular example, set $\cB$ is a vector space; we turn the reader's attention to \cite{LiLit05} for a survey of optimal policies in continuous-state MDPs.
\end{remark}

In general, we stipulate that the adversary's beliefs $B_t$ evolve according to some, potentially non-deterministic, memory-conscious learning mechanism 
\begin{equation}
\label{learnm}
\Prob(B_{t+1}=B)=f(s_0,\ldots,s_t,a_0,\ldots,a_t,B_0,t,B)\textrm{,}
\end{equation}
with $f:(\cup_{k=1}^\infty S^k\times A^k)\times\cB\times\{0,1,\ldots,T\}\times\cB\to\RR$, and where $B_0$ are the adversary's initial beliefs about the agent's reward. We note that equation \eqref{learnm} allows $B$ to be a general randomized estimator in the sense of \cite{LehCas98,Rob07}. While such a framework is sufficient for us to pose the problem of finding an optimal belief-induced policy, in practice we will require function $f$ to have a particular structure. Namely, if the belief updates are memoryless, i.e., $f$ depends only on $(s_t,B_t,a_t,B)$, the problem of finding an optimal belief-induced control policy is a problem of finding an optimal control policy on an MDP. It is given as follows:


\begin{problem}[Optimal deception]
\label{oripro}
Let $\cM=(S,A,P)$, $\cB$, and $L:S\times\cB\times A\to\RR$ be as above, and let $T\geq 0$, $s_0\in S$, $B_0\in\cB$. Find a control policy $\pi^*$ with such that \begin{equation}
\label{eqmdp}
\pi^*=\argmax_{\pi} \mathbb{E}\left[\sum_{t=0}^T L\left(s_t,B_t,\pi_t\right)\right]\textrm{,}
\end{equation}
subject to
\begin{equation}
\label{eff}
\begin{split}
\Prob(s_{t+1}=s) &=P(s_t,\pi_t,s)\textrm{,} \\
\Prob(B_{t+1}=B) &=f(s_t,B_t,\pi_t,B)\textrm{.}\end{split}
\end{equation}
\end{problem}

We note that, by \eqref{eqmdp}-\eqref{eff}, the optimal deception problem is a reward maximization problem on the MDP $\overline{\cM}=(S\times\cB,A,\overline{P})$, where $\overline{P}$ is given by \eqref{eff}. With such a model, it is well-known that the optimal policy $\pi^*$ is memoryless, i.e., $\pi^*_t$ depends solely on $s_t$, $B_t$, and $t$, and Problem \ref{oripro} is solvable by previously known and extensively discussed methods (see, e.g., \cite{Put05} for a detailed study).

\begin{remark}
The assumption that belief dynamics are memoryless, i.e., that beliefs change only by performing updates based on new information, and are not calculated directly from the history of information at every time step, holds for a wide variety of estimation techniques, most notably, online inverse reinforcement learning \cite{LiBur17}. In online inverse reinforcement learning, the parameter estimate $\hat{p}_t$ for a parametrized reward function $R_p$ is updated at every time step by setting $\hat{p}_{t+1}=\hat{p}_t+\alpha\nabla l(\hat{p}_t)$, where $l(p)$ is the estimated log-likelihood of the agent performing action $a_t$ at state $s_t$, if the reward function is given by $R_p$.
\end{remark}

Our running example, Example \ref{exa0}, can be simply stated in the MDP framework. We present such a setup in Section \ref{par}, where we also develop an optimal belief-induced policy for an agent. At the end of this section, let us, however, briefly discuss scenarios where agent's policy is constrained by predetermined specifications.

\subsection{Constraints on the Agent's Behavior}

As previously mentioned, the deception framework allows us to place additional constraints on the agent policies. Namely, Problem \ref{oripro} may be equivalently posed in the case where the set of all {\em permissible} policies $\pi$ is the proper subset of the set of all policies taking values in $A$. Notably, such a constraint may come from requirement that the agent follows a temporal logic specification; we direct the reader to a particularly detailed exposition given in \cite{Beletal17}. 

Clearly, a constraining specification may significantly lower the rewards that the agent is able to collect, making deception less effective. The extent to which a specification will make deception less effective is related to the extent to which it ``clashes'' with the agent behavior needed to successfully deceive the adversary. For instance, in our running example, a temporal logic specification may require that the agent never visits, or comes close to, one of the false candidate goals. Depending on the adversary's learning method, i.e., belief dynamics, such a specification may significantly limit the ability of the agent to impart an incorrect belief onto the adversary. On the other hand, a specification stating that the agent is never allowed to visit a particular state that is far away from the false candidate goals is not likely to significantly change the total accumulated reward gained by an optimal deceptive policy. Hence, the success of the deception does not depend on the agent visiting such a set.

We note that, while the problem of optimal deception with additional specifications is well-posed and formally given as an optimization problem analogous to \eqref{eqmdp}, with the additional constraint $\pi\in\cP$, such problems are generally computationally difficult to solve; we refer the reader to \cite{Dinetal11,Svoetal13} for recent work. 

The next section of the paper deals with settings where the agent lacks some knowledge about the adversary. In it, we do not emphasize the existence of possible constraints on the set of all policies $\pi$, with the understanding that the problem statements can easily be appended with such constraints. We return to such constraints in Section \ref{par}, where we illustrate the effectiveness of an optimally deceptive strategy with additional temporal logic specifications within the context of our running example.

\section{Lack of Knowledge on the Adversary in MDPs}
\label{lack2}

Problem \ref{oripro} poses the question of designing the optimal deceptive policy as a problem of finding an optimal policy on an MDP. While, as mentioned, such a problem can be solved by straightforward application of previously-known methods (e.g., value iteration \cite{Bel57}), its solution requires full knowledge of the dynamics and the reward function $L$ on the MDP $\overline{\cM}=(S\times\cB,A,\overline{P})$. Possession of such knowledge may not be realistic in some settings, especially considering that the motivating narrative for deception is of an adversarial scenario. As outlined in Section \ref{lack1}, we consider three categories in which the agent may lack knowledge on the belief-induced MDP $\overline{\cM}$:
\begin{enumerate}[1)]
\item knowledge of beliefs $B_t$ at any time $t$,
\item knowledge of the belief update mechanism \eqref{learnm}, and 
\item knowledge of the reward $L$ attained by the agent.
\end{enumerate}

Let us now consider each of the above possible knowledge limitations. For the sake of exposition, we treat these limitations separately, with the understanding that it is naturally possible that the knowledge about system $\cC_\cB$ is limited in more than one way at the same time, and that the above list of knowledge limitations is not exhaustive.

\subsection{Unknown Beliefs}
\label{subsube}

If the initial belief $B_0$ or current beliefs $B_t$, $t\geq 1$, are unknown, the system state $(s_t,B_t)$, evolving in $S\times\cB$, is partially observable. Namely, the agent knows $s_t$ at every time, but may not know $B_t$. Such a framework places the belief-induced system $\cC_\cB$ in the class of mixed-observability MDPs \cite{Ongetal09}, where the entirely observable part of the system state is $s_t$, and the entirely unobservable part is $B_t$. 

We will not describe the mixed-observability MDP framework in detail; we refer the reader to \cite{Ongetal09} for a more formal study. However, for the sake of phrasing the problem of designing an optimal deceptive policy, we note the fact that the agent cannot observe beliefs $B_t$ does {\em not} mean that the agent has no knowledge of $B_t$ whatsoever. Namely, if the agent possesses an initial probability distribution $Pr_0$ on possible $B_0\in\cB$, it may use the belief evolution \eqref{learnm} to obtain a probability distribution $Pr_1$ for $B_1\in\cB$, and by continuing onwards, distributions $Pr_t$ for $B_t$. 

The initial probability distribution $Pr_0$ depends on the agent's knowledge about the adversary's initial belief. If the agent has no knowledge about $B_0$, i.e., finds all beliefs in $\cB$ equally likely, the initial probability distribution $Pr_0$ is given by $Pr_0(B)=1/|\cB|$. On the other hand, if the agent knows that $B_0=B'$ for a particular $B'\in\cB$, $Pr_0$ is given by $Pr_0(B')=1$, $Pr_0(B)=0$ for all $B\neq B'$.

The problem of determining an optimally deceptive policy without belief observations is thus formalized as follows:

\begin{problem}[Optimal deception without belief observations]
\label{mompro}
Let $\cM=(S,A,P)$, $\cB$, $L:S\times\cB\times A\to\RR$, and a probability distribution $Pr_0:\cB\to[0,1]$ be as defined previously, and let $T\geq 0$, $s_0\in S$.

Find a control policy $\pi^*$, where $\pi^*_t=\pi^*(s_0,\ldots,s_t)$, such that \begin{equation}
\label{eqmomdp}
\pi^*=\argmax_{\pi} \mathbb{E}\left[\sum_{t=0}^T \sum_{B\in\cB} Pr_t(B)L\left(s_t,B,\pi_t\right)\right]\textrm{,}
\end{equation}
subject to
\begin{equation}
\label{eq2}
\begin{split}
\Prob(s_{t+1}=s) &=P(s_t,a,s)\textrm{,} \\
Pr_t(B) &=\sum_{B'\in\cB}Pr_t(B')f(s_t,B',\pi_t,B)\textrm{.}\end{split}
\end{equation}
\end{problem}

We emphasize that in the above problem, unlike in the other problems considered throughout this paper, $\pi^*$ needs to depend solely on the history of agent's positions since the beginning of the system run, as beliefs $B_0, \ldots, B_t$ are not known to the agent.

We note that \eqref{eqmomdp} produces a policy that maximizes the expected reward over the agent's {\em probability distribution} for the beliefs. Thus, if the agent's initial probability distribution $Pr_0$ is incorrect, i.e., does not satisfy $Pr_0(B_0)\approx 1$, the optimal policy in the sense of \eqref{eqmomdp} may not produce a high collected reward. On the other hand, if $Pr_0(B_0)=1$, i.e., the initial belief $B_0$ is known to the player, {\em and} the belief update mechanism \eqref{learnm} is deterministic, \eqref{eq2} will guarantee that $Pr_t(B_t)=1$ for all $t\geq 1$. Thus, the optimal policy in the sense of Problem \ref{mompro} will be the same as the optimal deceptive policy generated by the optimal deception problem (Problem \ref{oripro}).

Mixed-observability MDPs are a subclass of partially observable MDPs (POMDPs) and \cite{Ongetal09} provides an algorithm for determining an optimal policy for a mixed-observability MDP. Additionally, two surveys of algorithms for general POMDPs are given in \cite{Bra03,Mur00}. We omit further details of computing a policy that satisfies \eqref{eqmomdp}. We present an example of an approximation of such a policy in our running example, i.e., cops and deceptive robbers, in Section \ref{par}. In that example, the lack of knowledge on the adversary's beliefs naturally arises from the robbers not knowing the cops' estimate of their goal.

Finally, we note that a framework similar to Problem \ref{mompro} can be designed to deal with the more general case where $B_t$ is known by the agent for some $t$, but unknown for other times: if $B_t$ is observed, instead of evolving by \eqref{eq2}, $Pr_t$ is set to equal $1$ for the observed $B_t$, and $0$ for all other $B\in\cB$.

\subsection{Uncertain Belief Dynamics}

If the belief update mechanism $f$ in \eqref{learnm} is not entirely known, system governed by \eqref{eff} is transformed into an MDP with uncertain transition probabilities \cite{SatLav73}. In other words, it is known that \begin{equation}
\label{uncset}
f\in\cF=\{f^i:S\times\cB\times A\times\cB\to[0,1]~|~i\in I\}\textrm{,}
\end{equation}
where $I$ is an index set, and all $f^i$ satisfy $$\sum_{B'\in\cB}f'^i(s,B,a,B')=1\textrm{ for all }s\in S, B\in\cB, a\in A\textrm{.}$$
In such a case, the interest is to find a robust optimal policy, i.e., a policy that produces the best results for ``worst-case'' dynamics.

There are two basic variations of this problem: in one, the MDP transition probabilities, while uncertain, are the same at all times, while in the other, the transition probabilities are allowed to change over time, while remaining within the uncertainty set. Generally, the latter version is computationally easier to solve \cite{NilElG05}, and also allows for the possibility that the adversary does not learn in an entirely Markovian way, i.e., that the transition probabilities at time $t$ also depend on the entire history of the system states $(s_{t'},B_{t'})$ at times $t'=0,\ldots,t-1$, as long as the probabilities always remain within the uncertainty set. For this reason, in this paper we choose to describe this latter version of the problem.

\begin{problem}[Robust optimal deception with uncertain belief dynamics]
\label{uncdyn}
Let $\cM=(S,A,P)$, $\cB$, and $L$ be as before, and let $T\geq 0$, $s_0\in S$, $B_0\in\cB$. Let $\cF$ be as defined in \eqref{uncset}.

Find a control policy $\pi^*$ such that \begin{equation}
\label{eqfun}
\pi^*=\argmax_{\pi}\inf_{f_0,\ldots,f_{T-1}\in\cF} \mathbb{E}\left[\sum_{t=0}^T L\left(s_t,B_t,\pi_t\right)\right]\textrm{,}
\end{equation}
subject to dynamics
\begin{equation*}
\begin{split}
\Prob(s_{t+1}=s) &=P(s_t,a,s)\textrm{,} \\
\Prob(B_{t+1}=B) &=f_t(s_t,B_t,\pi_t,B)\textrm{.}\end{split}
\end{equation*}
\end{problem}

Problem \ref{uncdyn} describes an MDP with uncertain probabilities. Such a framework has been subject of substantial previous research and, for a wide variety of uncertainty sets, algorithms for efficiently computing the solution to \eqref{eqfun} have been proposed, based on robust dynamic programming. We turn the reader's attention to \cite{SatLav73,WhiEld94,NilElG05} for standard works. These works also deal with exploring the relationship of Problem \ref{uncdyn} to the problem where the transition probabilities are fixed prior to the system run. 

As in the case of unobservable beliefs, we present an example of an optimal policy for uncertain belief dynamics within our running example in Section \ref{par}. Such a setting is naturally motivated by the robbers not knowing the mechanism that the cops use to learn, i.e., update their beliefs.

\subsection{Uncertain Belief-Induced Reward}

If the knowledge of how the adversary's beliefs change the nominal reward $R$ into $L$ is not precise, system $\cC_\cB$ is transformed into an MDP with uncertain rewards: it is known that 
\begin{equation}
\label{uncl}
L\in\cL=\left\{L^i:S\times\cB\times A\to\RR~\vert~i\in I\right\}\textrm{,}
\end{equation}
where $I$ is an index set. For $s\in S$, $B\in\cB$, and $a\in A$, we denote $\cL(s,B,a)=\{L^i(s,B,a)~|~i\in I\}$. 

Analogously to the previous case, there are two basic cases: in the first, the rewards, while unknown to the agent, are fixed before the system run. In the second, they are allowed to be time-varying, while staying within the uncertainty set. The latter case is again computationally easier and is more easily motivated in our running example: the robbers' gains may differ every time they rob a bank. We present such a setup in Section \ref{par}. 

We formalize the latter version of the problem of deception with uncertain rewards as follows.

\begin{problem}[Robust optimal deception with uncertain rewards]
\label{uncrew}
Let $\cM=(S,A,P)$ and $\cB$ be as before, and let $T\geq 0$, $s_0\in S$, $B_0\in\cB$. For every $s\in S$, $B\in\cB$, $a\in A$, let $\cL$ be as defined in \eqref{uncl}, and assume that all $\cL(s,B,a)$ are bounded from below.

Find a control policy $\pi^*$ such that \begin{equation*}
\pi^*=\argmax_{\pi}\mathbb{E}\left[\sum_{t=0}^T \inf\cL\left(s_t,B_t,\pi_t\right)\right]\textrm{,}
\end{equation*}
subject to dynamics \eqref{eff}.
\end{problem}

The following result, with an obvious proof, reduces the problem of robust optimal deception with uncertain reward into a problem of finding an optimal policy in an MDP.

\begin{proposition}
Control policy $\pi^*$ is a solution of Problem \ref{uncrew} if and only if it is a solution of Problem \ref{oripro}, with $L$ in Problem \ref{oripro} replaced by $\inf\cL$. 
\end{proposition}

The other interpretation of uncertain rewards, where the reward is fixed a priori, is known as an imprecise-reward MDP \cite{DelMan07}. Robust optimal policies for imprecise-reward MDPs are usually based on solving a minimax regret optimization problem; we turn the reader's attention to \cite{DelMan07,RegBou11,XuMan09}, and references contained therein. While posing such a problem requires merely a small change from Problem \ref{uncrew}, its solution is generally computationally infeasible. Thus, we focus our attention in future sections to finding optimal policies for uncertain rewards in the sense of Problem \ref{uncrew}.

\section{Cops and Deceptive Robbers}
\label{par}

In this section, we provide a more thorough analysis of our running example, previously described in Figure \ref{fig0} and Example \ref{exa0}. We remind the reader that we consider a setup where the agent can collect a reward by reaching and remaining at a particular target state, while the adversary attempts to learn the location of the target state. If the adversary learns the location correctly, the agent's reward is removed, and it instead receives a significant penalty for arriving at its target state.

Effectively, the setup above is a version of the well-known {\em heaven and hell} example \cite{GefBon98} (we invite the reader to also see \cite{Gro06} for a detailed description of a version more similar to our setting), where there is no a priori hell, and heaven can become hell if the adversary finds out the heaven's location correctly. 

\begin{figure}[ht]
\centering
\begin{tikzpicture}[scale=0.7]
\draw[step=1cm] (0,0) grid (8,8);
\fill[red] (0,7) rectangle (1,8);

\fill[blue!50!black] (6,4) rectangle (5,5);
\fill[blue, opacity=0.3] (4,3) rectangle (5,4);
\fill[blue, opacity=0.3] (7,5) rectangle (6,6);
\node at (5.5,4.5) {\color{white} $TG$};
\node at (6.5,5.5) {$G_2$};
\node at (4.5,3.5) {$G_3$};
\end{tikzpicture}
\caption{An example of the gridworld in which the agent moves. Its true goal ($TG=G_1$) is marked in dark blue. The  agent knows the position of the true goal, while the adversary agent only knows that the goal is one of the three blue tiles $\{G_1,G_2,G_3\}$ on the map. The red tile is the agent's starting position.}
\label{fig1}
\end{figure}
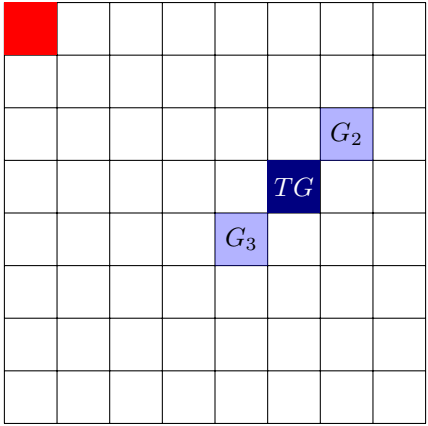

We assume that the agent moves in a gridworld $S$, shown in Figure \ref{fig1}. The actions available to the agent at any time are to go one tile north, south, east, west, or stay in place. (When the agent is at the edge of the grid state space, actions that would make it leave the state space are not available, or result in a prohibitively negative reward.) Given the agent's choice of action, the agent moves in the desired direction with probability $1$.

The agent's nominal reward $R$ is given by \begin{equation}
\label{nomrewhh}
R(s,a)=\begin{cases}
10 & \textrm{if } s=TG\textrm{,} \\
0 & \textrm{otherwise.}
\end{cases}
\end{equation}
Since the agent's goal is to maximize its reward over some time horizon, the nominal optimal behavior of the agent is to take the shortest path to $TG$, and then remain at $TG$ for the remaining time of the system run.

In the adversarial setting that we want to discuss, the adversary knows that the agent's reward function is given by \eqref{nomrewhh}. However, the adversary only knows of $k$ candidate tiles for $TG$, i.e., $G_1,\ldots,G_k$, where $TG\in\{G_1,\ldots,G_k\}$. Such a situation is illustrated in Figure \ref{fig1}, where $k=3$. We assume that the adversary has the knowledge of the agent's position and action at all times.

As proposed in Section \ref{defdec}, since the adversary is missing information about the exact location of the true goal, its belief space $\cB$ can be given by $\cB=\{1,\ldots,k\}$, where $B_t=i$ indicates that, at time $t$, the adversary believes that $TG=G_i$. The adversary uses the following memoryless mechanism for updating its beliefs: 
\begin{equation}
\label{upbehh}
\Prob(B_{t+1}=i|s_t,B_t,a_t)=m_1+m_2\textrm{,}
\end{equation}
with 
\begin{equation}
\label{upbehh1}
m_1=
\begin{cases}
0 & \textrm{if } B_t\neq i\textrm{,} \\
1-p & \textrm{if } B_t=i\textrm{,}
\end{cases}
\end{equation}
and
\begin{equation}
\label{upbehh2}
m_2=
\begin{cases}
0 & \textrm{if } d(s_{t+1},G_i)\geq d(s_t,G_i) \textrm{ and } s_{t+1}\neq G_i \textrm{ and } (B_t\neq i \textrm{ or } \#c_t\neq 0)\textrm{,}  \\
\frac{p}{\#c_t} & \textrm{if } d(s_{t+1},G_i)<d(s_t,G_i) \textrm{ or } s_{t+1}=G_i\textrm{,} \\
p & \textrm{if } B_t=i \textrm{ and } \#c_t=0\textrm{,}
\end{cases}
\end{equation}
where $d$ is the taxicab (i.e., $1$-norm) distance between two tiles, $\#c_t$ is the number of all $i\in\{1,\ldots,k\}$ such that $d(s_{t+1},G_i)<d(s_t,G_i)$ or $s_{t+1}=G_i$, and $p$ is a fixed parameter in $[0,1]$. We note that $s_{t+1}$ is a deterministic function of $(s_t,a_t)$, so the adversary's dynamics do not use any knowledge not available at the current time.

In plain words, \eqref{upbehh}-\eqref{upbehh2} state that the adversary's belief remains the same with probability $1-p$. The remaining $p$ are divided equally among all goal candidates, if any, which became closer to the agent as a result of the agent's last action. While such a learning mechanism is indeed simple, all of the work in this and the following section can be performed for any Markovian learning policy on a finite belief space. Additionally, even such a seemingly naive mechanism such as \eqref{upbehh}-\eqref{upbehh2} does guarantee that the adversary will eventually, with probability $1$, correctly learn the position of the true goal if the agent uses a nominal optimal control policy.

Finally, let us define the belief-induced reward $L$. It is modified from \eqref{nomrewhh} in such a way that, if the adversary's belief of the true position of the goal, the agent collects a negative reward:

\begin{equation}
\label{actrewhh}
L(s,B,a)=\begin{cases}
10 & \textrm{if } s=TG \textrm{ and } G_B\neq TG\textrm{,} \\
-10 & \textrm{if } s=TG \textrm{ and } G_B=TG\textrm{,} \\
0 & \textrm{otherwise.}
\end{cases}
\end{equation}

\subsection{Optimal Deceptive Policy}

As mentioned, belief update mechanism \eqref{upbehh}-\eqref{upbehh2} ensures that, if the agent is following the nominal optimal control policy $\pi^*_N$ (i.e., the policy that makes the agent take the shortest path to the true goal and remain there), the adversary will eventually correctly learn the position of the true goal, with probability $1$. Thus, $$\lim_{T\to +\infty}\mathbb{E}\left[\sum_{t=0}^T L(s_t,B_t,\pi^*_N(s_t))\right]=-\infty\textrm{,}$$
and, in particular, $$\lim_{T\to +\infty}\frac{\mathbb{E}\left[\sum_{t=0}^T L(s_t,B_t,\pi^*_N(s_t))\right]}{T}=-10\textrm{.}$$ Hence, not only is the nominal optimal control policy not optimal for the belief-induced system $\cC_\cB$, it is asymptotically the worst policy for such a system.

\begin{figure}[ht]
\centering
\includegraphics[width=0.48\textwidth]{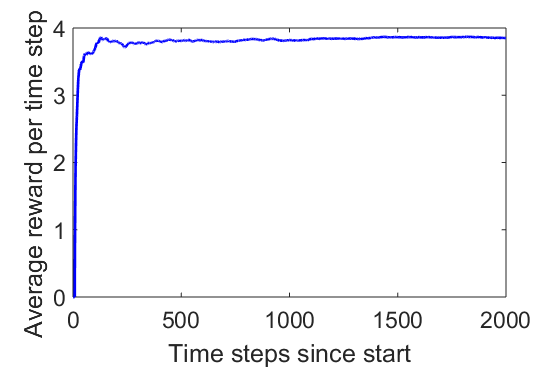}
\includegraphics[width=0.48\textwidth]{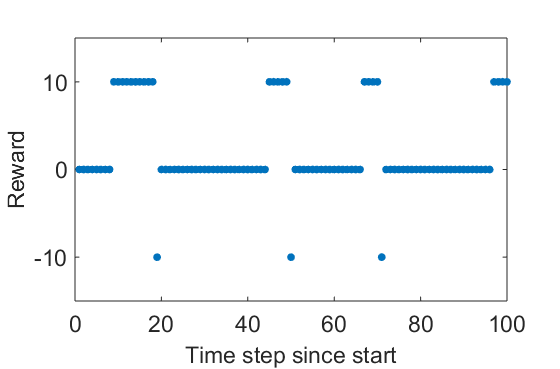}
\caption{The graph on the left side shows the accumulated reward at each time step, divided by number of time steps since start, obtained by the agent when using the optimal deceptive policy, over $100$ system runs. The graph on the right side shows the rewards obtained at each time step during the beginning of the system run.}
\label{fig2}
\end{figure}

It follows that there is clearly a need for determining the optimal belief-induced, or deceptive, policy $\pi^*_O$, which takes into account the adversary's beliefs $B_t$. As outlined in Section \ref{mdpsec}, $\pi^*_O$ is an optimal control policy for an MDP $\overline{\cM}=(S\times\cB,A,\overline{P})$ given by deterministic dynamics of the agent on $S$, dynamics \eqref{upbehh}-\eqref{upbehh2} on $\cB$, and reward function \eqref{actrewhh}. It can be constructed using any of the available algorithms for optimal control on MDPs (see \cite{Put05} for a detailed survey). The left side of Figure \ref{fig2} presents the average rewards \begin{equation}
\label{avgrew}
\sum_{t=0}^{T'} \frac{L(s_t,B_t,\pi^*_O(t))}{T'}
\end{equation}
obtained by the agent using such an optimal deceptive policy for $T'\leq T=2000$, with the adversary's belief change probability $p=0.1$, and the state space given as in Figure \ref{fig1}.

As the left side on Figure \ref{fig2} shows, unlike when using the nominal optimal control policy, the agent using the optimal deceptive policy $\pi^*_O$ collects, on average, a positive reward. We note that this average positive reward does not significantly depend on the initial state of the agent as the number of time steps grows large, and that the exact value of the average collected reward depends on the value of the belief change probability $p$. Nonetheless, the above simulation illustrates significant gains for the agent when using the optimal deceptive policy, rather than the nominal optimal policy which does not take the adversary's beliefs into account.

To help give some intuition to the agent's deceptive strategy, the right hand side of Figure \ref{fig2} shows the exact rewards that the agent collected during the first $100$ time steps in one system run. We note that the agent starts off by collecting a reward of $0$ for the first $8$ steps, until it reaches its goal $TG$. It then proceeds to remain at this goal until the adversary realizes that $TG$ is indeed the agent's true goal. After the adversary realizes the true goal and the agent collects a reward of $-10$, the agent leaves and tries to confuse the adversary by feigning that another one of the candidate goals is its goal, collecting a reward of $0$ during this period. Once the adversary is convinced incorrectly, the agent moves again to $TG$, and the process repeats. A video illustrating one typical system run is available at \url{https://bit.ly/2rtiygB}.

\subsection{Optimal Deception with Temporal Logic Specifications}

Let us briefly return to the setting where the agent is required to obey additional specifications while executing a deceptive policy. We give two examples. In the first one, the agent is not allowed to visit either of the two false goals, described by light blue tiles in Figure \ref{fig1}. In the second one, the agent is not allowed to visit the top false goal, but can visit the bottom.

\begin{figure}[ht]
\centering
\includegraphics[width=0.75\textwidth]{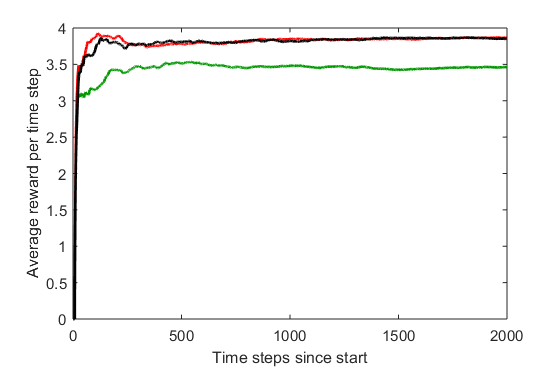}
\caption{The graphs show the accumulated reward at each time step, divided by number of time steps since start, obtained by the agent over $100$ system runs by three agent's policies. The graph in the green describes the rewards obtained by policy $\pi^*_{O1}$, where the agent is not allowed to visit either of the false goals. The graph in red describes the rewards obtained by $\pi^*_{O2}$, where the agent is not allowed to visit one of the false goals. The graph in the black describes the rewards obtained by $\pi^*_O$ --- an optimal deceptive policy without any specifications.}
\label{figadd}
\end{figure}

As described above, both of the above specifications yield constraints on the optimization problem \eqref{eqmdp} used in determining the optimal deceptive policy. In this case, the constraints are simple, and we can easily compute the optimal deceptive policies $\pi^*_{O1}$ and $\pi^*_{O2}$ for the two examples above, respectively. Figure \ref{figadd} shows the average rewards \eqref{avgrew} after each time step with each of the above strategies.

We note that the rewards obtained by policy $\pi^*_{O1}$ are significantly lower than those obtained by $\pi^*_{O2}$ or $\pi^*_O$. Such a difference arises from the fact that the specification underlying $\pi^*_{O1}$ is substantially interfering with the agent's ability to deceive the adversary --- not being able to go into any of the false goals makes it harder to convince the adversary that one of those goals is in fact the agent's objective. On the other hand, the rewards obtained by $\pi^*_{O2}$ are essentially the same as the rewards for $\pi^*_O$. While the specification underlying $\pi^*_{O2}$ restricts the agent's actions, it does not interfere with its ability to deceive the adversary --- if the agent is unable to visit one of the false goals, it will simply try to convince the adversary that the other false goal is the true objective. Thus, even though $\pi^*_{O2}$ may not be the same as $\pi^*_O$, the rewards that it obtains are not noticeably lower.

Let us now continue with the analysis of the optimal deceptive policies in the setting where the agent's knowledge about the adversary is not complete. In all further scenarios, we will not consider additional specifications placed on the agent's policies, but note that a similar analysis as above may be performed in all cases, and that the success of an agent at deceiving the adversary depends on interference between the specifications placed upon it and its optimal deceptive behavior.

\subsection{Optimal Deception with Imperfect Knowledge}
\label{opmik}

Having calculated and analyzed the optimal deceptive policy for an agent that is aware of all transition probabilities for the MDP $\overline{\cM}=(S,A,\overline{P})$, as well as reward $L$ and the adversary's beliefs $B_t$ at every time, let us now develop appropriate policies for each of the cases of imperfect knowledge discussed in Section \ref{lack2}.

\begin{figure}[ht]
\centering
\includegraphics[width=0.75\textwidth]{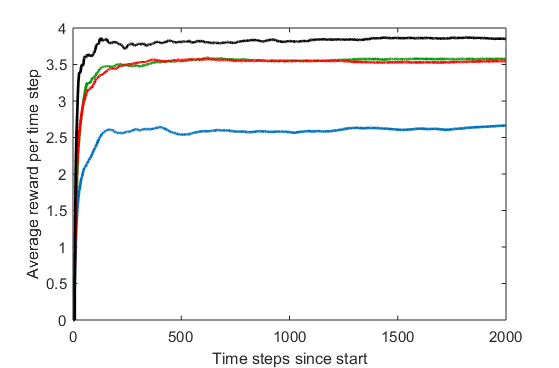}
\caption{The graphs show the accumulated reward at each time step, divided by number of time steps since start, obtained by the agent over $100$ system runs by four of agent's policies. The graph in light blue depicts the rewards obtained by an approximately optimal policy without belief observations. The graph in red depicts the rewards obtained by an optimal policy for the case of an uncertain learning parameter, while the graph in green depicts the rewards obtained by an optimal policy for the case where the collected rewards are a priori uncertain. The graph in black depicts the optimal deceptive policy with complete knowledge of the adversary.}
\label{figpom}
\end{figure}

In the case when the agent does not have any knowledge of the adversary's beliefs $B_t$, we showed in Section \ref{subsube} that an optimal deceptive policy is given by an optimal control policy of a mixed-observability MDP. As they depend on probability distributions on partially observed states, optimal policies for mixed-observability MDPs and POMDPs are generally difficult to compute exactly \cite{Chaetal12, Che88}. In this simulation, we used a randomized approximation of an optimal policy based on combining optimal actions for MDPs where beliefs are known, with weights corresponding to the probability distribution of the beliefs \cite{Caretal18}. The light blue graph in Figure \ref{figpom} describes average rewards \eqref{avgrew}, in analogy to the left side of Figure \ref{fig2}. As expected, the deceptive policy developed without belief observations performs worse than the optimal deceptive policy with perfect knowledge. However, such a policy still yields significantly better results than the nominal optimal policy.

\begin{remark}
\label{remlear}
We note that, if the agent is able to observe the reward that it collects at any given time, it might be able to use such information possible to deduce the adversary's beliefs. Namely, if the agent positioned at $TG$ collects a reward $L(TG,B_t,a)=-10$, it can immediately know that $G_{B_t}=TG$. On the other hand, if $L(TG,B_t,a)=10$, the agent can deduce that $G_{B_t}\neq TG$. In the simulation above, we do not consider that the agent performs such deductive reasoning. Nonetheless, a problem of designing an optimal strategy where such deduction is done is an interesting problem, as the agent has an incentive to go to those states where collecting a reward is more likely to inform it of the adversary's beliefs. We briefly discuss such a problem further in Section \ref{conc}.
\end{remark}

Let us now present the simulation results for the case of uncertain transition probabilities and uncertain rewards. In the former scenario, the agent does not know the true value of learning parameter $p$, and knows that it is between $0.05$ and $0.2$; it expects that $p$ may change at every time step, and designs the robust worst-case policy as a solution to Problem \ref{uncdyn}. In the simulation, $p$ is set to constant $0.1$, as before. In the latter scenario, the agent does not know the true reward that it will collect upon reaching $TG$ if the adversary did not learn its goal correctly, and believes it to be anywhere between $1$ and $20$. The agent expects that this reward may change every time it reaches $TG$, and designs the robust worst-case policy as a solution to Problem \ref{uncrew}. In the simulation, the collected reward is set to always equal $10$, as before.

We note that, while worst-case optimal policies, whose average collected rewards are illustrated in Figure \ref{figpom}, perform clearly worse than the optimal policy $\pi^*_O$ designed with complete information, the difference is less stark than for the case of unknown beliefs. The average collected reward \eqref{avgrew} for $\pi^*_O$ approaches around $3.9$, while the rewards for the optimal policies from Problem \ref{uncdyn} and Problem \ref{uncrew} approach around $3.5$. Such a property is a consequence of the simplicity of the reward function $L$: regardless of how quickly it believes the adversary is learning, or how large of a reward it may collect at goal, the agent has little motivation but to continue with the general behavior of reaching $TG$, waiting until the adversary learns of its goal, then moving away, confusing the adversary, and repeating the process. 

\begin{figure}[ht]
\centering
\includegraphics[width=0.75\textwidth]{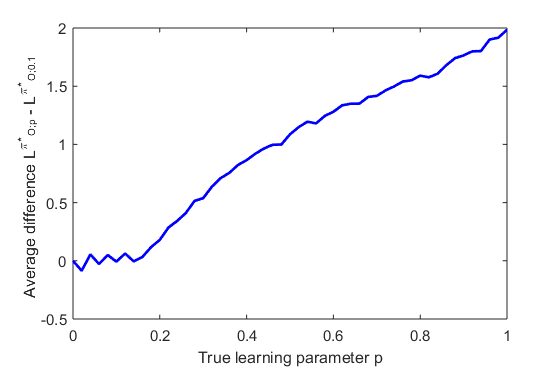}
\caption{The graph shows the difference between the average reward \eqref{avgrew} for $T'=T=2000$ collected by an agent when using the optimal policy $\pi^*_{O;p}$ for a learning parameter $p$ and the optimal policy $\pi^*_{O;0.1}$ for the learning parameter equaling $0.1$, in the scenario when the true learning parameter is $p$.}
\label{figdif}
\end{figure}

Finally, to illustrate the sensitivity of optimal deceptive policies under a small change in learning dynamics, we consider a scenario where the agent determines an optimal deceptive policy for $p=0.1$, with $p$ defined in \eqref{upbehh}-\eqref{upbehh2}, and uses it in the setting where the true learning parameter $p$ does not equal $0.1$. Figure \ref{figdif} shows the difference between the average reward $L^{\pi^*_{O;p}}$ obtained from \eqref{avgrew} by using the true optimal deceptive policy for a learning parameter $p$ and the average reward $L^{\pi^*_{O;0.1}}$ obtained by using the optimal deceptive policy for the case of $p=0.1$ in the setting where $p\neq 0.1$. While small imperfections in the agent's assumptions about true probabilities do not significantly influence the collected rewards (the differences can largely be attributed to randomness in instantiating each system run), the difference becomes clearer when $|p-0.1|>0.1$, and agent's incorrect assumption that $p=0.1$ can lead to significantly diminished rewards.

\section{Camouflage}
\label{camor}

In this section, we describe a scenario in which deception is attained through camouflage. In other words, the agent can, at a cost, choose to hide its location and thus influence the adversary's beliefs of its whereabouts. Camouflage techniques have been widely used in warfare --- for a historical survey of the use of camouflage in hostile actions, we direct the reader to \cite{For14}. In addition to motion camouflage \cite{Kwaetal15,RanIgl16,KimKim18}, we also note the emerging use of camouflage in disguising network traffic \cite{Guaetal01, WeiLep03}, as well as in hardware protection \cite{Cocetal14}.

The setting that we discuss in this section is the following: as in our running example, the agent moves deterministically along a gridworld $S$. Its nominal reward is given by a ``smoothened'' version of \eqref{nomrewhh}:
\begin{equation}
\label{nomrewca}
R(s,a)=\frac{10}{\|s-TG\|+1}\textrm{,}
\end{equation}
where $TG$ is again the location of the true goal.

This time, the nominal reward \eqref{nomrewca} {\em and} the position of $TG$ are both entirely known to the adversary. However, the adversary might not see the agent's current position. To formalize this setting, we stipulate that the agent has $10$ available actions: $A=\{\textrm{North}, \textrm{South}, \textrm{West}, \textrm{East}, \textrm{Stay in place}\}\times\{\textrm{Camo}, \textrm{No camo}\}$. We accordingly denote agent's actions by $a=(a^1,a^2)$. If the agent chooses to move without camouflage (``No camo''), its position becomes immediately visible to the adversary, but such an action does not incur any cost. If the agent uses camouflage, the adversary will, with a high probability $1-p$, not see the agent's movement, and hence its estimate of the agent's function will remain where it previously was. With probability $p$, the adversary will see the agent's position despite the camouflage. Hence, the adversary uses the following mechanism for updating its beliefs: 
\begin{equation}
\label{upbehhca}
\Prob(B_{t+1}=B|s_t,B_t,a_t)=
\begin{cases}
0 & \textrm{if } B\neq s_{t+1} \textrm{ and } a_t^2=\textrm{No camo,} \\
1 & \textrm{if } B=s_{t+1} \textrm{ and } a_t^2=\textrm{No camo,} \\
0 & \textrm{if } B\neq s_{t+1} \textrm{ and } B\neq B_t \textrm{ and } a_t^2=\textrm{Camo,} \\
p & \textrm{if } B=s_{t+1} \textrm{ and } B\neq B_t \textrm{ and } a_t^2=\textrm{Camo,} \\
1-p & \textrm{if } B=B_t \textrm{ and } B\neq s_{t+1} \textrm{ and } a_t^2=\textrm{Camo,} \\
1 & \textrm{if } B=B_t=s_{t+1} \textrm{ and } a_t^2=\textrm{Camo.} \\
\end{cases}
\end{equation}
We note that, naturally, the belief space $\cB$ equals $S$. 

The reward function $R$ is modified by the adversary's estimate of the agent's position: if the adversary believes that the agent is in position $B\in S$, then $L(s,B,a)=0$ for all $s\in S$ such that $d(s,B)\leq r$, $r>0$. However, while using camouflage may help the adversary, it incurs a cost. That is, the agent's collected reward is decreased by $c$, $c>0$, if $a^2=\textrm{Camo}$. Hence, the complete belief-induced reward is given by
\begin{equation}
\label{actrewca}
L(s,B,a)=\begin{cases}
-c & \textrm{if } \|s-B\|\leq r \textrm{ and } a^2=\textrm{Camo,} \\
0 & \textrm{if } \|s-B\|\leq r \textrm{ and } a^2=\textrm{No camo,} \\
R(s,a)-c & \textrm{if } \|s-B\|> r \textrm{ and } a^2=\textrm{Camo,}
\\
R(s,a) & \textrm{if } \|s-B\|> r \textrm{ and } a^2=\textrm{No camo.}
\end{cases}
\end{equation}
In our simulation, we use $r=1$ and $c=5$.

We note that the notion of camouflage, as an action that does not impact the agent's position in the state space $S$, but impacts solely the belief $\cB$, can be understood as a {\em gesture} in the sense of \cite{Bacetal16}. The framework of our paper can incorporate any such gestures that are made at the same time as actions that may lead to the agent's movement by making the action set a product of the set of all actions that may lead to the agent's movement and the set of all gestures.
\subsection{Optimal Deceptive Policy}

As in the previous example, the optimal nominal policy $\pi^*_N$ that solves \eqref{nomeqex} is for the agent is to take the shortest path to $TG$, without using camouflage, and then remain at $TG$ for the remaining time of the system run. It can be easily shown that $$\lim_{T\to +\infty}\frac{\mathbb{E}\left[\sum_{t=0}^T R(s_t,B_t,\pi^*_N(s_t))\right]}{t}=10\textrm{.}$$ Nonetheless, $\pi^*_N$ is not an optimal policy for the belief-induced scenario. Without using camouflage, the adversary will always be able to keep track of the agent, and thus eliminate its reward at every time step. In other words, $$L(s_t,B_t,\pi^*_N(s_t))=0$$ 
for all $t\geq 0$. On the other hand, if the agent knows that it is being watched, it can from time to time use its camouflage ability, even if such an action comes with a cost.

\begin{figure}[ht]
\centering
\includegraphics[width=0.55\textwidth]{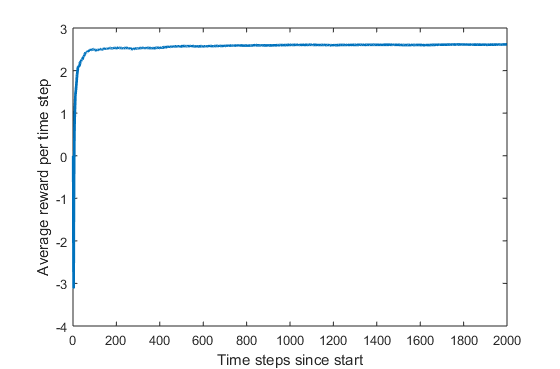}
\begin{tikzpicture}[scale=1]
\draw[step=1cm] (0,0) grid (5,5);
\fill[blue!50!black] (1,2) rectangle (2,3);
\draw[->, green, ultra thick] (1.5,2.5) -- (0.5,2.5);
\draw[->, green, ultra thick] (0.5,2.5) -- (0.5,1.5);
\draw[->, red, ultra thick] (0.5,1.5) -- (1.5,1.5);
\draw[->, red, ultra thick] (1.5,1.5) -- (1.5,2.5);
\draw[->, red, ultra thick] (1.5,2.5) arc(-135:215:0.5);
\node at (0,-0.5) {};
\end{tikzpicture}
\caption{The graph on the left describes the accumulated reward at each time step, divided by number of time steps since start, obtained by the agent when using the optimal deceptive policy, amalgamated from $100$ system runs. The graph on the right describes characteristic behavior of the agent upon reaching its true objective, when using the optimal deceptive policy. The objective $TG$ is depicted in blue. Movements performed without camouflage are denoted in green, while camouflaged movements are denoted in red.}
\label{fig3}
\end{figure}
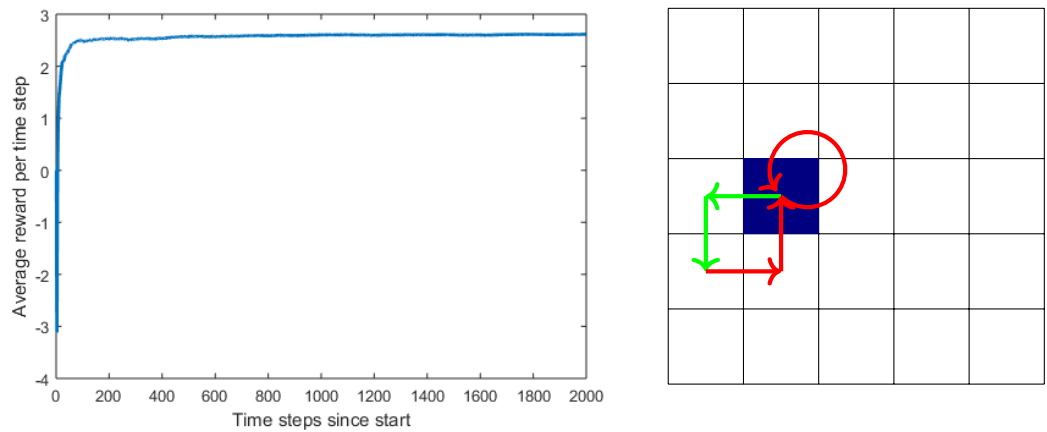

Let us thus compute an optimal deceptive policy $\pi^*_O$. The gridworld in this example is a $5\times 5$ grid shown on the right hand side of Figure \ref{fig3}. We note that the belief-induced state space $S_\cB=S\times\cB$ then has $|S|^2=625$ elements. The left hand side of Figure \ref{fig3} presents the average rewards \eqref{avgrew} collected by the agent using the optimal policy $\pi^*_O$, with the camouflage strength parameter $p$ set to $p=0.1$. As seen in Figure \ref{fig3}, the optimal deceptive policy $\pi^*_O$ again, as expected, significantly outperforms $\pi^*_N$.  To illustrate the intuition behind the optimal agent's policy, let us briefly describe its typical behavior. While approaching its goal, the agent uses camouflage in an effort to deceive the adversary about its position and prevent it from reducing the collected reward. Once it reaches the objective, the agent continues to use camouflage while remaining at the same position, until the adversary learns its true location, after which the agent leaves the objective without using camouflage, thus allowing the adversary to follow its movement. Once it is far enough from the objective, depending on the adversary's vision radius $r$, the agent then turns back towards the objective and moves back to it under camouflage. Such behavior is illustrated on the right hand side of Figure \ref{fig3}, and a typical system run is shown in a video available at \url{https://bit.ly/2K1B1HV}.

We note that, as in the cops and deceptive robbers example, the adversary is simple to trick. As it is memoryless, the adversary does not remember that the agent keeps repeating the same sequence of actions: reaching the goal, remaining at it until discovered, then moving away to confuse the adversary, and ultimately returning to the goal. One could certainly design a ``smarter'' adversary to deal with both of the presented examples, and we briefly discuss such design in the final section of the paper. Nonetheless, we emphasize that adversaries that we presented in the above two examples will, with probability $1$, correctly learn the agent's intentions if the agent uses a nominal optimal policy. Additionally, the optimal deceptive policies used by the agent firmly correspond to the natural strategy that one would use, or attempt to use, when trying to deceive an adversary. Thus, the policies presented above verify that the formal meaning of deception as defined in Section \ref{defdec} corresponds to common intuition.

\section{Conclusions and Future Work}
\label{conc}

This paper provided a formal definition of deception and deceptive strategies in optimal control, applicable to a wide variety of adversarial scenarios. The definition that we presented rests on the introduction of the belief space of an adversary who is trying to learn agent's intentions, as well as on encoding the agent's control objective and adversary's influence on the objective through nominal and belief-induced rewards, respectively. The problem of devising an optimal policy for a belief-induced reward, i.e., an optimal deceptive policy, results in an optimal control problem on a state space that is a product of the agent's original state space and the adversary's belief space. In the context of this paper, we primarily focused on an agent whose behavior is governed by a discrete-time discrete-state MDP. Assuming that the adversary's learning process is memoryless, i.e., its beliefs at every time are updated based solely on the current belief and the agent's current action, the problem of optimal design of a deceptive policy is an optimal control problem in an MDP. In the natural case when not all information about the adversary's learning process is known, this problem converts into a question of control policy design for partially observable or uncertain MDPs. Through two examples in which deception naturally arises --- cops and robbers scenario and a setting in which an agent may use camouflage --- we showed that an optimal deceptive policy, as defined in the paper, corresponds to common intuition behind deceptive behavior, and achieves significantly better results for the deceiving agent than if deception was not used.

Let us briefly discuss avenues of future work; this paper is merely a first step in analyzing deception within an optimal control setting. Firstly, while the above MDP setting encompasses a wide variety of scenarios, two of which were discussed within this paper, there is clearly interest in analysis of deceptive strategies in other settings. Namely, the adversaries' learning mechanisms in the context of the examples of this paper were simple --- an adversary that operates without memory and has only finitely many potential beliefs can be easily deceived by a single ``trick'' performed time and time again during the system run (e.g., the ``leaving the objective'' strategy illustrated on the right hand side of Figure \ref{fig3}). Determining deceptive strategies against adversaries with more advanced learning mechanisms would have significant use. A particularly interesting mechanism is the online inverse reinforcement learning \cite{LiBur17} --- an algorithm particularly designed to discover the agent's reward function. Secondly, the ability to analyze deception and design deceptive strategies in a less constrained setting than that of MDPs would enable us to place deception in a wide variety of applications within the framework developed in this paper. As mentioned in the introduction, deception is used in a variety of contexts, and framing those instances of deception --- currently described in vastly differently fashions --- within a common structure would be a significant breakthrough.

Within the context of future work, we want to emphasize again that the three settings of lack of agent's knowledge that we discussed in this paper are in no way exhaustive. In particular, all three settings are essentially static. As mentioned in Remark \ref{remlear}, the agent is not discovering anything new about the adversary during the system run; it either knows a particular element of the adversary throughout the run, or it does not know it at any time. Such a stipulation might be too constraining for some scenarios. For instance, even if the agent does not usually know the adversary's belief, it can partly infer it from the collected rewards: as mentioned in the cops and robbers example, if the agent collects a reward of $-10$ when at $TG$, it can directly deduce the adversary's belief, even if it does not directly observe it. Allowing the agent to  learn about the adversary during the system run would yield a more realistic set of scenarios to be handled by the presented framework, and could add another dimension to the agent's policy. Namely, in addition to trying to achieve its nominal objective and deceive the opponent, the agent would now be motivated to direct its actions in such a way that it learns about the adversary. Such a research direction would result in an immediate connection between deception and reinforcement learning, thus opening the door for discussion of a complex, but perhaps more realistic, strategy in adversarial long-term scenarios: first learn about the adversary, and then use the learned information to deceive it.

In addition to computationally determining optimal deceptive policies for a variety of settings with full or limited information, designing a way to describe agent behavior that {\em makes deception succeed} --- the salient features that all optimal deceptive policies for a particular scenario need to have in order to collect maximal belief-induced rewards --- would be of significant interest. For instance, in both of the examples presented in this paper, optimal deceptive policies relied on the agent moving away from the objective once its intention has been uncovered, then acting in a way that aims to trick the adversary about its intentions, finally followed by returning to its objective once the adversary has been tricked. However, for each example there were multiple policies that exhibited such behavior. Hence, while finding an optimal deceptive policy, or enumerating all optimal deceptive policies, may be sufficient in order to deceive an adversary, analyzing and learning from deceptive behavior in one scenario in order to determine deceptive behavior in a similar scenario would require us to describe the set of deceptive policies in understandable terms. Such a question broadly falls within the research effort on explainable artificial intelligence \cite{Swaetal91,vanetal04,Coretal06}.

Finally, we note that the definition of deceptive policies presented in this paper rests on the assumption that the state space is endowed with a reward function that encodes the agent's objective. It is naturally possible to consider a different class of objectives. Namely, objectives encoded in temporal logic specifications are particularly intuitive for a variety of applications. In the current paper, we considered temporal logic specifications solely as an addition to the reward objective: if such specifications exist, they serve to reduce the set of agent's admissible policies. Such a framework is significantly different from the setup in which the agent's sole objective is given by a temporal logic specification, which the adversary then attempts to learn. It is possible to encode temporal logic specifications as reward functions in a product MDP \cite{Woletal12}, where the maximal expected accumulated future reward of an agent positioned at a particular state corresponds to the probability of satisfying the temporal logic specification. Additionally, our framework already permits the belief set $\cB$ to consider any properties of interest to the adversary, including temporal logic specifications. However, the notion of belief-induced rewards does not neatly carry over to the temporal logic framework: since the rewards in the product MDP considered in \cite{Woletal12} are produced by the MDP transition probabilities, changing those rewards based on adversary's beliefs would imply that the adversary is able to change the agent's dynamics. In contrast, in the framework of this paper, agent's dynamics are invariant to adversary's beliefs. Thus, it would be meaningful to devise a new natural notion of deception for scenarios with objectives expressed in temporal logic, instead of assimilating it into the rewards framework by the process considered in \cite{Woletal12}.

\section*{Acknowledgement}
The authors thank Steven Carr for coding and running the simulation of an deceptive strategy in the scenario where the observer beliefs are unknown, in Section \ref{opmik}. This work was funded by grants W911NF-16-1-0001 from the Defense Advanced Research Projects Agency, FA8650-15-C-2546 from the Air Force Research Laboratory, and W911NF-15-1-0592 from the Army Research Office.

\bibliographystyle{IEEEtran}
\bibliography{refs}
\end{document}